\documentclass[11pt]{amsproc}
\usepackage{amsfonts}
\usepackage{amsmath,amstext,amsbsy,amssymb}
\newtheorem{theorem}{Theorem}[section]
\newtheorem{lemma}[theorem]{Lemma}

\newtheorem{corollary}[theorem]{Corollary}

\theoremstyle{definition}

\theoremstyle{remark}

\numberwithin{equation}{section}

\newcommand{\la}{\lambda}

\begin{document}

\title[On vertex operators and Jack functions]{On vertex operator realizations of
Jack functions}
\author{Wuxing Cai}
\address{School of Sciences,
South China University of Technology, Guangzhou 510640, China}
\email{caiwx@scut.edu.cn}
\author{Naihuan Jing$^*$}
\address{Department of Mathematics,
   North Carolina State University,
   Raleigh, NC 27695-8205, USA}
\email{jing@math.ncsu.edu}
\thanks{*Corresponding author}
\thanks{Jing gratefully acknowledges the support from
NSFC's Overseas Distinguished Youth Grant (10801094).}

\begin{abstract}
On the vertex operator algebra associated with rank one lattice we
derive a general formula for
 products of vertex operators in terms of generalized homogeneous symmetric functions.
 As an application we realize Jack symmetric functions of
 rectangular shapes as well as marked rectangular shapes.
\end{abstract}

\maketitle

\vskip 0.1in

\section{Introduction}

Classical symmetric functions play important roles in various
areas of mathematics and physics, and they admit many different
formulation. Starting with Bernstein's work \cite{Ze}, vertex
operators have been used in constructing several families of
symmetric functions such as Schur and Schur's Q-functions
\cite{J1} as well as Hall-Littlewood symmetric functions
\cite{J2}. Though one can still define certain family of vertex
operators associated with Macdonald polynomials similar to the
Schur case, the products of the vertex operators are in general no
longer equal to the Macdonald polynomials. At least in the case of
two rows the transition function from the basis of generalized
homogeneous symmetric functions (product of one-row Mcdonald
polynomials) involves with hypergeometric series of type
$_4\phi_3$. Except for a few cases \cite{J3, Z}, it has been an
open problem to find the transition function from the generalized
homogeneous functions to Macdonald symmetric functions \cite{G,
J3}.

 On the other hand, in the vertex representations of affine Lie algebras,
Lepowsky and Wilson \cite{LW} have long posted the important
problem on whether certain products in the representation space
are linearly independent. In the special cases of level three
representations, this problem can be solved using Rogers-Romanujan
identities (see also \cite{LP} for the homogeneous case). Later it
was realized \cite{J3} that the vertex operators at level three
are actually related with vertex operators associated to certain
Jack polynomials \cite{J3}, namesly, the half vertex operators are
actually the generating function of the one-row Jack functions.
Thus it is also an interesting question to study the linear
independence problem for those vertex operators associated to Jack
functions.

Motivated by \cite{FF} we define a new type of vertex operators
associated with Jack functions in the vertex operator algebra of
rank one lattice. We will call them Jack vertex operator since the
product of identical modes of this vertex operator will be shown
to be Jack functions of rectangular shapes. It is interesting that
the contraction functions for products of the new vertex operators
are of the form $\prod_{i<j} (z_i-z_j)^{2\alpha}$ instead of
$\prod_{i<j} (z_i-z_j)^{\alpha}$ for the Jack parameter
$\alpha^{-1}$ (which are for $Y_1(z)$ in section 3.2), as expected
from experiences with Schur and Hall-Littlewood cases. It turns
out that one really needs this new form of vertex operators
(vertex operator $X(z)$ in section 3.2) to generate rectangular
Jack functions. At the special case of Schur functions
($\alpha=1$), our new vertex operators provide another formula for
the rectangular shapes.

We also study the problem of linear independence for the new
vertex operators in the case of Jack functions. We show that under
certain conditions the set of vertex operator products are indeed
a basis for the representation space (see \cite{FF} for a similar
statement). We achieve this by deriving a Jacobi-Trudi like
formula for the Jack vertex operators, and then we reprove
Mimachi-Yamada's theorem \cite{MY} that the product of the vertex
operators are Jack functions for the rectangular shapes, and then
we further generalize this formula to the case of marked
rectangular shapes, i.e., rectangular shapes minus a row of boxes
at the lower left corner. This general case includes \cite{JJ} as
special cases.

This paper is organized as follows. In section 2 we recall some
necessary notions of symmetric functions. In section 3 we first
review the vertex operator approach based on the second author's
work on Schur functions, then we define the Jack vertex operators
and give an explicit formula of the vertex operator products and a
Jacobi-Trudi like formula in terms of tableaux. In section 4 we
provide a detailed analysis of certain matrix coefficients of
vertex operators and prove the theorem of realizing Jack functions
of rectangular and marked rectangular shapes.

\section{Jack functions}

 We recall some basic notions about
symmetric functions following the standard reference \cite{M}.
   A partition $\lambda$ is a sequence
   $\lambda=(\lambda_1,\lambda_2,\cdots,\lambda_s)$ of nonnegative
   integer such that
   $\lambda_1\geq\lambda_2\geq\cdots\geq\lambda_s$; the set of all partitions is denoted by $\mathcal {P}$;
   we sometimes write $\lambda$ as
    $\lambda=(1^{m_1}~2^{m_2}\cdots)$,
   where $m_i$ is the multiplicity of $i$ occurring in the sequence.
   The number of non-zero $\lambda_i$'s is called the length of
   $\lambda$, denoted by $l(\lambda)$, and the weight $|\lambda|$
    is defined as
$\lambda_1+\cdots+\lambda_s$. We also recall that the dominance
order is defined by
   comparing the partial sums of the parts. For two partitions
   $\la$ and $\mu$ of the same weight, if
    $\lambda_1+\cdots+\lambda_i\geq\mu_1+\cdots+\mu_i$ for all
   $i$, one says that $\lambda$ is greater than $\mu$ and denoted as
   $\lambda\geq\mu$; conventionally, $\lambda>\mu$ means $\lambda\geq\mu$ but
   $\lambda\neq\mu$.
   For
   $\la=(\la_1,\la_2,\cdots)=(1^{m_1},2^{m_2},\cdots)$, $\mu=(\mu_1,\mu_2,\cdots)=(1^{n_1},2^{n_2},\cdots)$
the notation $\la-\mu$ means $(\la_1-\mu_1,\la_2-\mu_2,\cdots)$,
$\mu\subset'\la$ means that $n_1\leq m_1,n_2\leq m_2,\cdots$, and
$\la\backslash\mu$ denotes the partition
$(1^{m_1-n_1}2^{m_2-n_2}\cdots)$. We also define
$\binom{m(\la)}{m(\mu)}=\binom{m_1}{n_1}\binom{m_2}{n_2}\cdots$, and
$\lambda\cup\mu=(1^{m_1+n_1}~2^{m_2+n_2}\cdots)$.

The ring $\Lambda$ of symmetric functions over $\mathbb Z$ has
various linear $\mathbb Z$-bases indexed by partitions: the monomial
symmetric functions $m_{\la}=\sum x^{\la_1}_{i_1}\cdots
x_{i_k}^{\la_k}$, the elementary symmetric functions
$e_{\la}=e_{\la_1}\cdots e_{\la_k}$ with $e_n=m_{(1^n)}$, and the
Schur symmetric functions $s_{\la}$. The power sum symmetric
functions $p_{\la}=p_{\la_1}\cdots p_{\la_k}$ form a $\Bbb Q$-basis.

Let $F=\mathbb{Q}(\alpha)$ be the field of rational functions in
indeterminate $\alpha$. The Jack polynomial is a special orthogonal
symmetric function under the following inner product. For two
partitions $\la, \mu \in \mathcal P$ the scalar product on
$\Lambda_F$ is given by
\begin{align} \label{def}
<p_{\la}, p_{\mu}>=\delta_{\la, \mu}\alpha^{-l(\la)}z_\lambda
\end{align}
where $z_\lambda=\prod_{i\geq 1}i^{m_i}m_i!$, $m_i$ is the
occurrence of integer $i$ in the partition $\la$ , and $\delta$ is
the Kronecker symbol. Here our parameter $\alpha$ is chosen as the
reciprocal to the usual convention in view of our vertex operator
realization.

In \cite{M} Macdonald proved the existence and uniqueness of what is
called the Jack function as a distinguished family of orthogonal
symmetric functions $P_{\lambda}(\alpha^{-1})$ with respect to the
scalar product (\ref{def}) in the following sense:
$$
P_{\la}(\alpha^{-1})=\sum_{\la\geq\mu}c_{\la
\mu}(\alpha^{-1})m_{\mu}
$$
in which $c_{\la \mu}(\alpha^{-1})\in F$ , $\la, \mu \in \mathcal
P$, and $c_{\la \la}(\alpha^{-1})=1$. Let
$Q_{\lambda}(\alpha^{-1})=b_{\la}(\alpha^{-1})P_{\lambda}(\alpha^{-1})$
be the dual Jack function.

It is known that the special case $Q_{(n)}(\alpha^{-1})$, simplified
as $Q_n(\alpha^{-1})$, can be written explicitly:
\begin{equation}\label{E:homogeneous1}
Q_n(\alpha^{-1})=\sum_{\lambda\vdash
n}\alpha^{l(\lambda)}z_\lambda^{-1}p_{\lambda}.
\end{equation}
For a partition $\la$, we will denote
$q_{\la}(\alpha^{-1})=Q_{\la_1}(\alpha^{-1})Q_{\la_2}(\alpha^{-1})\cdots
Q_{\la_l}(\alpha^{-1})$.

According to Stanley \cite{S} the $q_\la$'s also form another basis
of $\Lambda_F$, and they are dual to that of $m_{\la}$. Hence the
transition matrix from $Q_\la$'s to $q_\la$'s is the transpose of
that from $m_\la$'s to $P_\la$'s. Explicitly, we have

\begin{lemma} \label{L:triangular}
For any partition $\lambda$, one has
$$q_{\lambda}(\alpha^{-1})=\sum_{\mu\geq\lambda}c'_{\lambda\mu}Q_{\mu}(\alpha^{-1})$$
$$Q_{\lambda}(\alpha^{-1})=\sum_{\mu\geq\lambda}d'_{\lambda\mu}q_{\mu}(\alpha^{-1})$$
where $d'_{\lambda,\mu}\in F$,with $c'_{\lambda\mu}=c_{\mu\lambda}$
and  $c'_{\lambda\lambda}=d'_{\lambda\lambda}=1$.
\end{lemma}

\section{Vertex operators and symmetric functions}

Vertex operators can be used to realize several classical types of
symmetric functions such as Schur and Hall-Littlewood polynomials
\cite{J1, J2}. There are some partial progress  towards realizations
of Macdonald polynomials \cite{J1, Z}. In order to discuss the Jack
case, we will use the standard vertex algebra technique and recall
the construction of lattice vertex operator algebra for rank one
case.

\subsection{Representation space $V$ and transformation to $\Lambda_ C$}
For a positive integer ${\alpha}$, the complex Heisenberg algebra
$H_{{\alpha}}=\displaystyle\bigoplus_{n\neq 0}$ $\mathbb{C}$
$h_n+\mathbb{C}$ $c$ is the infinite dimensional Lie algebra
generated by $h_n$ and $c$ subject to the following defining
relations:
   $$[h_m,h_n]=\delta_{m,-n}{\alpha}^{-1}mc,~~[h_m,c]=0.$$
We remark that the integer ${\alpha}$ is included for identification
with Jack inner product. If it is clear from the context, we will
omit the subscript ${\alpha}$ in $H_{{\alpha}}$ and simply refer it
as $H$.

  It is well known that $H$ has a unique canonical
  representation given as follows.  The representation space can be
  realized as the infinite dimensional vector space
  $V_0=Sym(h_{-1},h_{-2},\cdots)$, the symmetric algebra over
    $\mathbb{C}$, generated by $h_{-1},h_{-2},\cdots$.
The action of $H$ is given by
\begin{eqnarray*}
      & &h_n.v= {\alpha}^{-1}n\frac{\partial}{\partial h_{-n}}v\\
      & &h_{-n}.v=h_{-n}v\\
      & &c.v=v
\end{eqnarray*}

% \textbf{4.1.The degree on V.}\vskip 0.2in

To simplify the indices we enlarge the space $V_0$ by the group
algebra of $\mathbb Z$. Let $V=V_0\otimes$ $\mathbb{C}[\mathbb{Z}]$,
where
  $\mathbb{C}[\mathbb{Z}]$ is the group algebra of $\frac{1}{2}\mathbb Z$ with
  generators
  $\{e^{nh}|n\in\frac{1}{2}\mathbb{Z}\}$. We define the action of the
  group algebra as usual with the multiplication given by
  $e^{mh}e^{nh}=e^{(m+n)h}$. We also define the action of
  $\partial=\partial_h$ on $\mathbb C[Z]$ by
  $\partial.e^{mh}=me^{mh}$.
The space $V_0$ is $\mathbb Z$-graded. The enlarged space $V$ is
  doubly $\mathbb Z$-graded as follows. Let $\la=(\lambda_1,\lambda_2,\cdots)$ be a partition, for
$v=h_{-\la}\otimes e^{nh}\in
  V$,
  define the degree of $v$ by $nd(v)=(|\la|,n)$, where we have used the usual notation
$h_{-\lambda}=h_{-\lambda_1}h_{-\lambda_2}\cdots$. For convenience,
  we consider the degree of zero element to
  be of any value.\vskip 0.1in
  The vertex operator space $V$ has a canonical scalar product. For
any polynomials $P, Q$ in the $h_i$'s we have

 \begin{eqnarray*}
  & &\langle h_n P,Q\rangle=\langle P,h_{-n}Q\rangle\\
  & &\langle1,1\rangle=1\\
  & &\langle e^{mh},e^{nh}\rangle=\delta_{m,n}
  \end{eqnarray*}
Thus, for partitions $\la,\mu$ we have $$\langle h_{-\lambda}
\otimes e^{mh},h_{-\mu} \otimes
  e^{nh}\rangle=z_\lambda {\alpha}^{-l(\lambda)} \delta_{\lambda,\mu} \delta_{m,n}.$$
  We define a linear map $T$: $V=\sum_{s\in \frac{1}{2}\mathbb
Z}V_s\mapsto\Lambda_\mathbb{Q}$ by:
$$T:h_{-\lambda}\otimes e^{sh} \mapsto p_\lambda.$$
We remark that the restriction of $T$ on $V_s=V_0\otimes e^{sh}$ is
a bijection preserving the products.
\subsection{Jack vertex operator on $V$ and an explicit formula }

For a complex parameter $a$ we let the vertex operator $Y_{a}(z)$
acts on V via the
  generating series:
  $$Y_{a}(z)=exp\Big(\sum_{n=1}^\infty \frac{z^n}{n}{\alpha}h_{-n}\Big)
  exp\Big(\sum_{n=1}^\infty \frac{z^{-n}}{-n}a{\alpha}h_n\Big)=\sum_n Y_a(n)z^{-n}.$$
  For $\alpha\in \frac12\mathbb Z$, we define $X(z)=Y_2(z)exp(2{\alpha}lnz\partial_h+h)$, i.e.
  $$X(z)=exp\Big(\sum_{n=1}^\infty \frac{z^n}{n}{\alpha}h_{-n}\Big) exp(2{\alpha}lnz\partial_h+h)
  exp\Big(\sum_{n=1}^\infty \frac{z^{-n}}{-n}2{\alpha}h_n\Big)$$
  where the middle term acts as follows:
\begin{equation}
exp(2{\alpha}lnz\partial_h+h).e^{sh}=z^{(s+\frac{1}{2})2{\alpha}}e^{(s+1)h}.
\end{equation}
  The operator $X_n$ on V is defined as the component of $X(z)$:
  $$X(z)=\sum_{n\in\mathbb{Z}} X_nz^{-n}.$$

%$X(z)$ is the main operator we will study about.  But we will use
%vertex operators without middle term to generate a basis of
 %$\Lambda_{F}$.
 For simplicity we consider a special case of the vertex operator
 $Y_a(z)$, and let %$\Lambda_{\mathbb{C}}$\vskip 0.1in
\begin{align*}
Y(z)&=Y_0(z)=exp\Big(\sum_{n=1}^\infty
\frac{z^n}{n}{\alpha}h_{-n}\Big)=\sum
Y_{-n}z^n,\\
Y^*(z)&=Y_0^*(z)=Y(z^{-1})^*=exp\Big(\sum_{n\geq
1}\frac{z^{-n}}{n}{\alpha}h_n\Big)=\sum Y_n^*z^n,
\end{align*}
where we took the dual of $Y(z)$. We remark that one can also use
the operator $Y_{a}(z),(a\neq 0)$ in place of $Y(z)$, and most
proofs will remain the same.\vskip 0.0in

We note that when $\alpha=1$ the vertex operator $X(z)$ differs from
the Schur vertex operator \cite{J1} or its truncated form is not
Bernstein operator for Schur functions.

 To simplify
the notations, for partition $\lambda=(\lambda_1,\cdots,\lambda_s)$,
we
  denote the product
  $X_{-\lambda_1}\cdots X_{-\lambda_s}$ simply as
  $X_{-\lambda}$, and similarly for $Y_{-\la}$.

We will first give an explicit formula for the vertex operator
products. For this purpose, we need the following

\begin{lemma} \label{L:action}%{Lemma 4.4.1}
For the the creation part and annihilation part of $X(z)$ we have
\begin{equation} \exp\Big(\sum_{n=1}^\infty
\frac{z^n}{n}{\alpha}h_{-n}\Big)=\sum_{\lambda\subset\mathcal
{P}}z_{\lambda}^{-1}{\alpha}^{l(\lambda)}h_{-\lambda}z^{|\lambda|}\end{equation}
 and

\begin{equation}exp\Big(\sum_{n=1}^\infty
\frac{z^{-n}}{-n}2{\alpha}h_n\Big).h_{-\lambda}=\sum_{\mu\subset\la}\binom{m(\lambda)}
{m(\mu)}(-2)^{l(\mu)}z^{-|\mu|}h_{-\lambda\backslash\mu}\end{equation}
for any partition $\lambda$,
\end{lemma}
Proof: The first one is a direct computation:
$$exp\Big(\sum_{n=1}^\infty
\frac{z^n}{n}{\alpha}h_{-n}\Big)=\prod_{n\geq
1}exp(\frac{z^n}{n}{\alpha}h_{-n})=\prod_{n\geq
1}\sum_{i\geq0}\frac{1}{i!}\frac{z^{ni}}{n^i}{\alpha}^ih_{-n}^i=\sum_{\lambda\in\mathcal
{P}}z_{\lambda}^{-1}{\alpha}^{l(\lambda)}h_{-\lambda}z^{|\lambda|}.$$
For the second one, let $h_n^{(i)}={h_n^i}/{i!}$. The Heisenberg
canonical commutation relation implies that
\begin{equation}
h_n^{(i)}.h_{-n}^m=(\frac{n}{{\alpha}})^i\binom{m}{i}h_{-n}^{m-i}.
\end{equation}

Using this we have
\begin{align*}
exp(
\frac{z^{-n}}{-n}2{\alpha}h_n).h_{-n}^m&=\sum_{i\geq0}(\frac{{2{\alpha}z}^{-n}}{-n})^ih_n^{(i)}.h_{-n}^m\\
& =\sum_{i\geq0}z^{-ni}(-2)^i\binom mi h_{-n}^{m-i}.
\end{align*}
For $\lambda=(1^{m_1}2^{m_2}\cdots$), we have

$$exp\Big(\sum_{n=1}^\infty
\frac{z^{-n}}{-n}2{\alpha}h_n\Big).h_{-\lambda}=\prod_{n\geq 1}exp(
\frac{z^{-n}}{-n}2{\alpha}h_n).h_{-1}^{m_1}h_{-2}^{m_2}\cdots$$$$=\prod_{n\geq
1}\Big(\sum_{i_n\geq0}z^{-ni_n}(-2)^{i_n}\binom{m_n}{i_n}h_{-n}^{m_n-i_n}\Big)
=\sum_{\mu\in\mathcal{P}}\binom{m(\lambda)}{m(\mu)}
(-2)^{l(\mu)}z^{-|\mu|}h_{-\lambda\backslash\mu},$$ where the sum
runs through all partitions $\mu\subset'\la$.

\medskip

In the following $\underline{\la}=(\la^1,\la^2,\cdots,\la^s)$
denotes that $\underline{\la}$ is a sequence of partitions
$\la^1,\la^2,\cdots,\la^s$.
 \begin{theorem} \label{T:formula1}For
integer $s\geq 1$, we have
\begin{align} \nonumber
&X_{-\la_s}\cdots X_{-\la_1}.e^{nh}\\ \label{E:explicitformula}&=
\sum_{\underline{\mu},\underline{\nu}}\prod_{i=1}^s
\frac{(-2{\alpha})^{l(\nu^i)}}{z_{\nu^i}}
\binom{m(\mu^{i-1})}{m(\mu^i \backslash
\nu^i)}\frac{h_{-\mu^s}}{(-2)^{l(\mu^s)}}\otimes e^{(n+s)h},
\end{align}
where the sum is over $ \underline{\mu}=(\mu^1,\mu^2,\cdots,\mu^s)$
and $\underline{\nu}=(\nu^1, \nu^2,\cdots,\nu^s)$ such that
$\nu^i\subset'\mu^i$, $|\mu^i|=|\la^{(i)}|-i(2n+i){\alpha}$,
$\mu^i\backslash\nu^i\subset'\mu^{i-1}$, $\mu^0=(0)$, here
$\la^{(i)}=(\la_1,\cdots,\la_i)$ is a subpartition of $\la$.
\end{theorem}

Proof: We use induction on $s$. Note that $\nu^1=\mu^1$, it is
trivial for the case of $s=1$ . Applying the annihilation part and
the middle term of $X(z)$ to Eq. (\ref{E:explicitformula}), we find
the right side
\begin{align*}
&\sum_{\mu\subset'\mu^s}\sum_{\underline{\mu},\underline{\nu}}\prod_{i=1}^s
\frac{(-2{\alpha})^{l(\nu^i)}}{z_{\nu^i}}
\binom{m(\mu^{i-1})}{m(\mu^i \backslash
\nu^i)}\cdot \\
&\qquad\qquad\frac{(-2)^{l(\mu)}}{(-2)^{l(\mu^s)}}\binom{m(\mu^{s})}{m(\mu)}
\frac{z^{(2n+2s+1){\alpha}}}{z^{|\mu|}}h_{-\mu^s\backslash\mu}
e^{(n+s+1)h}.
\end{align*}

Replacing $\mu$ with $\mu_s\backslash\mu$, it becomes
\begin{align*}
&\sum_{\mu\subset'\mu^s}\sum_{\underline{\mu},\underline{\nu}}\prod_{i=1}^s
\frac{(-2{\alpha})^{l(\nu^i)}}{z_{\nu^i}}
\binom{m(\mu^{i-1})}{m(\mu^i \backslash
\nu^i)}\cdot\\
&\qquad\qquad (-2)^{-l(\mu)}\binom{m(\mu^{s})}{m(\mu)}
\frac{z^{(2n+2s+1){\alpha}}}{z^{|\mu^s\backslash\mu|}}h_{-\mu}
e^{(n+s+1)h}.
\end{align*}
Thus we have,
\begin{align*}
&X_{-\la_{s+1}}X_{-\la_s}\cdots X_{-\la_1}.e^{nh}
=\sum_\nu\sum_{\mu\subset'\mu^s}\sum_{\underline{\mu},\underline{\nu}}\prod_{i=1}^s
\frac{(-2{\alpha})^{l(\nu^i)}}{z_{\nu^i}}
\binom{m(\mu^{i-1})}{m(\mu^i \backslash
\nu^i)}\cdot\\
&\qquad\qquad\qquad\qquad\qquad\qquad(-2)^{-l(\mu)}\binom{m(\mu^{s})}{m(\mu)}
\frac{{\alpha}^{l(\nu)}}{z_{\nu}}h_{-\mu\cup\nu} e^{(n+s+1)h},
\end{align*} where in the first sum $\nu$ is restricted to:
$|\nu|-|\mu^s\backslash\mu|=\la_{s+1}-(2n+2s+1){\alpha}$. Changing
variables $\nu=\nu^{s+1}$,~~$\mu=\mu^{s+1}\backslash\nu^{s+1}$, we
finish the proof.

\medskip

Note that a necessary condition for $X_{-\la_s}\cdots
X_{-\la_1}.e^{nh}\neq 0$ is that
$k_i=|\la^{(i)}|-i(2n+i){\alpha}\geq 0$ for $i=1,2,\cdots,s$.

\subsection{Generalized Jacobi-Trudi's theorem}

The generalized homogeneous symmetric function with parameter
$\alpha$ can be written as follows. For $n\in\mathbb{Z}$,
\begin{equation}\label{E:homogeneous2}
H_n(\alpha^{-1})=\sum_{\lambda\vdash
n}\alpha^{(\lambda)}z_\lambda^{-1}h_{-\lambda}.
\end{equation}
The image of it under $T$ is $Q_n(\alpha)$ by \ref{E:homogeneous1}.
Clearly one has $Y_a(-n).1=H_n({\alpha}^{-1})$. Notice that by
definition $H_n({\alpha}^{-1})=0$ for $n<0$. Combining these,we have
the following statement after a simple computation (see also Lemma
(\ref{L:triangular}))

%\textbf{4.2. The normalization} of $X(z_1)\cdots X(z_s)$.\vskip
%0.1in

\begin{lemma}  \label{L:homogeneous2} %{Lemma 4.3.2 }
For any positive integer $s$, we have
\begin{equation*}
Y(z_1)\cdots Y(z_s)e^{mh}= \sum_{n_1\geq0, \cdots,
n_s\geq0}H_{n_1}({\alpha}^{-1})\cdots
H_{n_s}({\alpha}^{-1})z_1^{n_1}\cdots z_s^{n_s}.
\end{equation*}
In particular, $Y_{-r_1}\cdots Y_{-r_s}\cdot
1=H_{r_1}({\alpha}^{-1})\cdots H_{r_s}({\alpha}^{-1})$ for any
integers $r_1, \cdots, r_s$. Moreover, for partition $\la$
$$T(Y_{-\la}\cdot 1)=q_{\la}({\alpha}^{-1}).$$
\end{lemma}\vskip 0.3in

Thus the vectors $\{Y_{-\lambda}e^{mh}\}$ ($\lambda\in\mathcal P$
and $m\in\frac{1}{2}\mathbb Z$) form a linear basis in the
representation space, corresponding to the basis of generalized
homogeneous polynomials.

To proceed further we define the normalization of vertex operators,
which helps to separate the singular part. The normalization of
$X(z_1)\cdots X(z_s)$ is defined as
\begin{align*} &:X(z_1)\cdots
X(z_s):\\
&=exp\Big(\sum_{n=1}^\infty
\frac{z_1^n+\cdots+z_s^n}{n}{\alpha}h_{-n}\Big)
exp\Big(\sum_{n=1}^\infty
\frac{z_1^{-n}+\cdots+z_s^{-n}}{-n}2{\alpha}h_n\Big)A,
\end{align*}
where
$A=e^{sh}exp({2\alpha}(lnz_1+\cdots+lnz_s)\partial_h)(z_1\cdots
z_s)^{\alpha}$
%$A=exp({2\alpha}lnz_1\partial_h+h)\cdots
%exp({2\alpha}lnz_s\partial_h+h)$.
%We remark that for convenience our
%normal order excludes the middle terms.
Similarly when the normal
product is taken on mixed product of $X(z)$, $Y(z)$, and $Y^*(z)$,
one always moves the annihilation operators to the right.

 We define the lowering operator $D_{i}$ on the bases
of symmetric functions by
$$D_i(H_{\la})=H_{\la_1}\cdots H_{\la_{i-1}}H_{\la_i-1}H_{\la_{i+1}}\cdots
H_{\la_l}.$$ The rising operator is defined by
$R_{i,j}=D_i^{-1}D_j$. Like the raising operator $R_{ij}$  the
lowering operator $D_i$ is not always invertible, one needs to make
sure that each application of $D_i$ is non-zero.

\begin{lemma} \label{L:product-of-X}
For $s\geq1$, we have
\begin{eqnarray*}
X_{-\lambda_1}X_{-\lambda_2}\cdots X_{-\lambda_s}.
e^{nh}&=&\Big(\prod_{1\leq i<j\leq s}(D_i-D_j)^{2{\alpha}}\cdot
\prod_{1\leq i\leq
s}D_i^{(n+\frac{1}{2})2{\alpha}}\Big)\\
&\cdot& H_{\lambda_1}({\alpha}^{-1})\cdots
H_{\lambda_s}({\alpha}^{-1})\otimes e^{(n+s)h}.
\end{eqnarray*}
\end{lemma}
Proof: Observe that
\begin{eqnarray*}
&&exp\Big(\sum_{n=1}^\infty \frac{z^{-n}}{-n}2{\alpha}h_n\Big)
exp\Big(\sum_{n=1}^\infty \frac{w^n}{n}{\alpha}h_{-n}\Big)\\
%&=&exp\Big(\sum_{n=1}^\infty \frac{w^n}{n}{\alpha}h_{-n}\Big)
%exp\Big(\sum_{n=1}^\infty \frac{z^{-n}}{-n}2{\alpha}h_n\Big)
%exp\Big(\sum_{n=1}^\infty-\frac{(w/z)^n}{n}2{\alpha} \Big)\\
&=&exp\Big(\sum_{n=1}^\infty \frac{w^n}{n}{\alpha}h_{-n}\Big)
exp\Big(\sum_{n=1}^\infty \frac{z^{-n}}{-n}2{\alpha}h_n\Big)
\Big(1-\frac{w}{z}\Big)^{2{\alpha}}~~(|w|<|z|)  % ~~~(4.1)
\end{eqnarray*}

Using induction on $s$ and the normalization
 we have \vskip 0.1in
\begin{equation*}
X(z_1)\cdots X(z_s)=\prod_{1\leq i<j\leq
s}(z_i-z_j)^{2{\alpha}}:X(z_1)\cdots X(z_s):,  %(4.2)
\end{equation*}

Applying the action on $e^{nh}$, we have

\begin{align*} \nonumber
X(z_1)\cdots X(z_s).e^{nh}&=\prod_{1\leq i<j\leq
s}(z_i-z_j)^{2{\alpha}}\prod_{i=1}^s z_i^{(n+\frac{1}{2})2{\alpha}}\\
&\cdot\prod_{i=1}^s \Big(\sum_{n\geq
0}H_n({\alpha}^{-1})z_i^n\Big)\otimes e^{(n+s)h},
\end{align*}
where $(z_i-z_j)^{2{\alpha}}=z_i^{2{\alpha}}(1-z_j/z_i)^{2{\alpha}}$
if there is an expansion. Taking the coefficient of $z^{\lambda}$ we
obtain the statement.%\vskip 0.3in

\medskip
In concern with the operator in Lemma \ref {L:product-of-X}, we have
the following result on the square of Vandermonde determinant.

\begin{lemma} \label{L:square-of-Vandermonde}
Let $V(X_s)=V(x_1,\cdots,x_s)=\prod_{1\leq i<j\leq s}(x_i-x_j)$,
$s\geq 2$. For $V(X_s)^2$, the coefficient of the term
$\prod_{i=1}^sx_i^{s-1}$ is $(-1)^{s(s-1)/2}s!$, and the coefficient
of $x_kx_s^{-1}\prod_{i=1}^sx_i^{s-1}$ $(k=1,\cdots,s-1)$ is
$-(-1)^{s(s-1)/2}(s-1)!$.
\end{lemma}

Proof: The Vandermonde determinant $V(X_s)$ is the determinant of
$M=(x_j^{i-1})_{s\times s}$. Then
$V(X_s)^2=det(MM^T)=det(p_{i+j-2})_{s\times s}$, where
$p_n=x_1^n+\cdots+x_s^n$. The product $\prod_{1\leq i\leq
s}x_i^{s-1}$ only appears in the (sub-diagonal) term
$p_{s-1}p_{s-1}\cdots p_{s-1}=p_{s-1}^s$ of the determinant, thus
the coefficient is $(-1)^{s(s-1)/2}s!$. Similarly the term
$x_kx_s^{-1}\prod_{1\leq i\leq s}x_i^{s-1}$ ( $k=1,\cdots,s-1$) only
appears in the term of the form $p_{s-2}p_sp_{s-1}^{s-2}$ of the
determinant, so the coefficient is $-(-1)^{s(s-1)/2}(s-1)(s-2)!$.
   %\vskip 0.2in
\medskip

For any partition $\lambda$ and a fixed parameter $\alpha$, we set
\begin{equation}
H_{\lambda}(\alpha)=H_{\la_1}(\alpha)\cdots H_{\la_l}(\alpha).
\end{equation}
Clearly the set of vectors $H_{\la}({\alpha}^{-1})e^{mh}$ forms an
$F$-basis of the vertex operator space $V$. Under the map $T$, the
vector $H_{\la}({\alpha}^{-1})$ is the symmetric function
$q_{\la}({\alpha}^{-1})$. For fixed ${\alpha}\in\mathbb N$ and
$m\in\mathbb Z$, we define $\mathcal P_{{\alpha}, m}$ to be the
set of partitions $\lambda$ such that $\lambda_i-\lambda_{i+1}\geq
{\alpha}$ and $\lambda_l\geq \frac 12(2m+1)\alpha$. The following
result is a generalization of Jacobi-Trudi theorem for our vertex
operator basis.

\begin{theorem} \label{T:vertex-basis}
The set of products $X_{-\la}e^{mh}$ ($\la\in\mathcal P_{2\alpha,
m}, m\in\mathbb Z$) forms an $\mathbb F$-basis in the vertex
algebra $V$. Moreover one has, for a partition $\lambda$ of length
$l$ and $\lambda_l\geq (2m+1)\alpha$,
\begin{equation}
X_{-\lambda}
e^{mh}=\sum_{\mu\geq\la}a_{\la\mu}({\alpha}^{-1})H_{\mu-(2m+1){\alpha}{\bf
1}-2{\alpha}\delta}e^{(m+l(\la))h},
\end{equation}
where $\mu$ runs through the compositions such that
$a_{\la\la}({\alpha}^{-1})=1$, ${\bf 1}=(1, \cdots, 1)\in \mathbb
N^l$ and $\delta=(l-1, l-2, \cdots, 1, 0)$.
\end{theorem}
Proof. For any partition $\la$ of length $l$, we can rewrite Lemma
\ref{L:product-of-X} in terms of raising operators.
\begin{align} \nonumber
X_{-\lambda_1}X_{-\lambda_2}\cdots X_{-\lambda_l}\cdot
e^{mh}&=\Big(\prod_{1\leq i<j\leq l}(1-R_{ij})^{2{\alpha}}\cdot
\prod_{1\leq i\leq
l}D_i^{(m+l-i+\frac{1}{2})2{\alpha}}\Big)\\
&\cdot H_{\lambda_1}({\alpha}^{-1})\cdots
H_{\lambda_l}({\alpha}^{-1})\otimes
e^{(m+l)h}. % ~~~~~(4.4)
\end{align}
The raising operators map $H_{\la}$ into $H_{\mu}$ with $\mu\geq
\la$, and the product
$$\prod_{1\leq i<j\leq
l}(1-R_{ij})^{2{\alpha}}=1+\sum_{e\neq 0}
\pm\prod_{i<j}R_{ij}^{e_{ij}},$$ where $e_{ij}$ are non-negative
exponents.
 The equality is clear now. When $\lambda_l\geq
(2m+1)\alpha$, the composition $\lambda-(2m+1){\alpha}{\bf
1}-2{\alpha}\delta$ is a partition.  Then when $e\neq 0$, all the
terms in the sum differ from $H_{\la-(2m+1){\alpha}{\bf
1}-2{\alpha}\delta}e^{(m+l(\la))h}$, which shows that transition
matrix from the basis $H_{\mu}e^{nh}$ to the set $X_{-\la}e^{mh}$ is
triangular and has ones on the diagonal. On the other hand, any
vector $X_{-\la}e^{mh}$ can be expressed as a linear combination of
$X_{-\mu}e^{(m-l(\lambda))h}$, where $\mu\geq
\lambda+(2m-2l(\lambda)+1){\alpha}{\bf 1}+2{\alpha}\delta$. Hence
the set forms a basis of the vertex operator algebra. \vskip 0.3in

We will see that in certain cases the vectors $X_{-\la}e^{mh}$ are
actually Jack symmetric functions.

\subsection{Jack functions of rectangular shapes}

We observe that

\begin{lemma}\label{L:weight} For $\lambda\in\mathcal {P}$,
$n\in\frac{1}{2}\mathbb{Z}$, set $u=X_{-\lambda}.e^{nh}$, then
$nd(u)=(|\lambda|,n+l(\lambda))$, if and only if
$n=-\frac{l(\lambda)}{2}$.
\end{lemma}

Proof. First for $v$ such that $nd(v)=(m,n)$, by definition we have
$$nd(X_{-k}.v)=(m+k-(n+\frac{1}{2})2{\alpha},n+1)$$ And then we have
$nd(X_{-\lambda}.v)=(m+|\lambda|-{\alpha}(2n+l(\lambda))l(\lambda),n+l(\lambda))$.
The result follows.\vskip 0.3in

\begin{lemma} \label{L:partition} %{Lemma 4.2.2}
Let $\lambda=\Big((k+1)^s,(k)^t\Big)$ be a partition with $t\in
\mathbb{Z}$$_{>0}$, $k,s\in \mathbb{Z}$$_{\geq 0}$. Then for any
partition $\mu$ satisfying $|\mu|=|\lambda|$ and $l(\mu)\leq
l(\lambda)$, we have $\mu \geq \lambda.$
\end{lemma}
Proof: Let $\lambda=(\lambda_1,\cdots,\lambda_{s+t})$. If $s=0$ it
is obviously true. So we can assume that $s\geq 1$. Let
$\mu=(\mu_1,\cdots,\mu_{s+t})$ be another partition with
$\mu_1\geq \mu_2\geq\cdots \geq \mu_{s+t}\geq 0$. The assumption
says that $\mu_1\geq k+1$. If we do not have $\mu \geq \lambda$,
there should be an $r$, $1\leq r<s+t$ such that
$$\sum_{i=1}^r \mu_i \geq \sum_{i=1}^r\lambda_i$$ and
$$\sum_{i=1}^{r+1} \mu_i<\sum_{i=1}^{r+1}\lambda_i.$$ It then follows
that $k\geq \mu_{r+2},\mu_{r+3},\cdots$, but $k\leq
\lambda_{r+2},\lambda_{r+3},\cdots$. Subsequently
$$\sum_{i=1}^{s+t} \mu_i<\sum_{i=1}^{s+t}\lambda_i ,$$ a
contradiction with $|\mu|=|\lambda|$. \vskip 0.3in
 Next we consider
the mixed products. The following result is an easy computation by
vertex operator calculus.
\begin{lemma}\label{L:OPE}%\label{P:OPE} {Proposition 4.3.5}
The operator product expansion of mixed product is given by
\begin{align*}
&Y^*(w_1)\cdots Y^*(w_t)X(z_1)\cdots X(z_s)= %exp(lnz_s\partial_h+h)\cdots(lnz_1\partial_h+h)\cdot
\\
&:-: \prod_{1\leq i<j\leq
s}(z_i-z_j)^{2{\alpha}}\prod_{j=1}^s\prod_{i=1}^t(1-z_jw_i)^{-{\alpha}}
,
\end{align*}
where $:-:=:X(z_1)\cdots X(z_s)Y(w_1)\cdots Y(w_t):$.
\end{lemma}

 Now we can prove the main theorem.
\begin{theorem} \label{L:squareJack}  For partition $\lambda=((k+1)^s,(k)^t)$ with $k\in
\mathbb{Z}_{\geq 0}$, $s\in \mathbb{Z}_{> 0}$, $t\in \{0,1\}$, we
have
$$T(X_{-\lambda}e^{-(s+t)h/2})=c(\alpha)Q_{\lambda}({\alpha}^{-1}),$$
where $c(\alpha)$ is a rational function of $\alpha$, and
$c(1)=(-1)^{s(s-1+2t)/2}s!$.
\end{theorem}
We remark that when $\la$ is a rectangular tableau (i.e. $t=0$) the
result was first proved by Mimachi and Yamada
\cite{MY} using
differential operators. When $s=t=1$, it was proved in \cite{JJ}.
Another important phenomenon is that when $\alpha=1$, we obtain a
new vertex operator formula for the rectangular shapes and marked
rectangular shapes.

Proof. For $\lambda=((k+1)^sk^t)$, let
$u=T(X_{-\lambda}.e^{-(s+t)h/2})$. Note that $u$ is a linear
combination of $Q_\mu(\alpha^{-1})$'s with $\mu\geq\la$, by Lemma
\ref{L:product-of-X} and Lemma \ref{L:triangular}. By Lemma
\ref{L:weight}, we need to show that $u$ is orthogonal to
$Q_{\mu}({\alpha}^{-1})$, for all $\mu$ such that $\mu\vdash
kt+s(k+1)$ and $\mu\neq \lambda$. This will be done in the
following. As for the coefficient, we have applied Lemmas
\ref{L:square-of-Vandermonde} and \ref{L:triangular} as well as
Lemma \ref{L:product-of-X}, which confirms that $u$ is non-zero and
the coefficient $c(\alpha)$ satisfies the given formula.\vskip 0.0in

To prove the orthogonality, consider the two cases of $\mu$:\vskip
0.0in (1) By Corollary \ref{L:partition} and Lemma
\ref{L:triangular}, for $\mu<\lambda$, or $\mu$ is incomparable with
$\lambda$, it follows easily that $u$ is orthogonal to
$Q_{\mu}({\alpha}^{-1})$.

(2) If partition $\mu=(m_1,\cdots,m_r)$ and $\mu>\lambda$ (or $\mu$
are incomparable with $\lambda$), then it follows that $m_1>k+1$. By
Lemma \ref{L:triangular} and Lemma \ref{L:homogeneous2}, to prove
that $u$ is orthogonal to $Q_{\mu}({\alpha}^{-1})$, we just need to
prove the following product is zero:

$$\langle ((X_{-(k+1)})^s(X_{-k})^t.e^{-(s+t)h/2},Y_{-m_r}\cdots
Y_{-m_1}.e^{(s+t)h/2}\rangle$$
$$=\langle Y^*_{-m_1}\cdots
Y^*_{-m_t}(X_{-(k+1)})^s(X_{-k})^t,e^{(s+t)h/2}\rangle, $$ which
equals to the coefficient of $w_1^{-m_1}\cdots w_t^{-m_t}(z_1\cdots
z_s)^{k+1}(z_{s+1}\cdots z_{s+t})^k$ in the following expression

\begin{align} \nonumber
&\langle Y^*(w_1)\cdots Y^*(w_r)X(z_1)\cdots
X(z_{s+t}).e^{-(s+t)h/2},e^{(s+t)h/2}\rangle\\ \nonumber
&=z_{s+t}^{(-\frac{s+t}{2}+\frac{1}{2})2{\alpha}}z_{s+t-1}^{(-\frac{s+t}{2}+1+\frac{1}{2})2{\alpha}}\cdots
z_1^{(-\frac{s+t}{2}+s+t-1+\frac{1}{2})2{\alpha}}\cdot \\
\nonumber&\qquad\qquad\qquad \prod_{s+t\geq j>i\geq
1}(1-z_jz_i^{-1})^{2{\alpha}}\prod_{i=1}^{s+t}(1-z_jw_i^{-1})^{-{\alpha}}\\
\label{E:innerproduct} &=\pm \prod_{1\leq i\neq j\leq
s+t}(1-z_jz_i^{-1})^{{\alpha}}\prod_{j=1}^{s+t}\prod_{i=1}^r(1-z_jw_i^{-1})^{-{\alpha}},
\end{align}
where we have used Lemma \ref{L:OPE}. This coefficient in Eq.
(\ref{E:innerproduct}) is zero by Lemma \ref{L:contraction}, which
we will prove in the next section. Hence the theorem is proved.

In general we have the following result.

\begin{corollary}  For partition $\lambda=((k+1)^s,(k)^t)$ with $k\in
\mathbb{Z}_{\geq 0}$, $s\in \mathbb{Z}_{> 0}$, $t\in \{0,1\}$, and
$r\in\mathbb{Z}$, $\alpha\in\mathbb{Z}_{>0}$ such that
$r\alpha\leq k+\delta_{t,0}$ we have
$$T(X_{-\lambda}e^{-\frac{s+t-r}{2}h})=cQ_{\lambda-r\alpha {\bold 1}}(\alpha^{-1})$$
where $c$ is a nonzero constant, and $\bold
1=(1,1,\cdots,1)\in\mathbb{Z}^{s+t}$.
\end{corollary}

Proof: The proof is essentially the same as that of Theorem
\ref{L:squareJack}. The condition $r\alpha\leq k+\delta_{t,0}$ is
included to make sure that
$T(X_{-\lambda}e^{-\frac{s+t-r}{2}h})\neq 0$ (see Remark of
Theorem \ref{T:formula1}).

\begin{lemma}\label{L:contraction} %{Lemma 4.3.6}
The contraction function $H_{{\alpha}}(Z_s, W_t)$ does not contain
terms like $w_1^{-m_1}\cdots w_t^{-m_t}z_1^{k_1}\cdots z_s^{k_s}$ if
$m_1>k_i (i=1,2,\cdots,s)$ where
$$H_{{\alpha}}(Z_s, W_t)=\prod_{s\geq i\neq j\geq 1}(1-z_jz_i^{-1})^{{\alpha}}\prod_{j=1}^s
\prod_{i=1}^t(1-z_jw_i^{-1})^{-{\alpha}}.$$
\end{lemma}

\section{Analysis of $H_{{\alpha}}(Z_s, W_t)$}
We have the following lemma to {\it split} $H_{{\alpha}}(Z_s,W_t)$:
\vskip 0.1in
  \begin{lemma}\label{L:split}For positive integers $r,s$, $i\neq j$, there are
  non-negative integers $f_m$ and $g_n$ such that:
$$\Big(\frac{1-z_iz_j^{-1}}{1-z_iw^{-1}}\Big)^r \Big(\frac{1-z_jz_i^{-1}}{1-z_jw^{-1}}\Big)^s
  =\sum_{m=1}^r\Big(\frac{1-z_iz_j^{-1}}{1-z_iw^{-1}}\Big)^mf_m+\sum_{n=1}^s
  \Big(\frac{1-z_jz_i^{-1}}{1-z_jw^{-1}}\Big)^ng_n.$$

\end{lemma}
Proof: For simplicity we denote
$a=\frac{1-z_iz_j^{-1}}{1-z_iw^{-1}}$,
$b=\frac{1-z_jz_i^{-1}}{1-z_jw^{-1}}$
 , it can be verified directly that $ab=a+b$. Repeatedly using
this, we can write $a^rb^s$  into the wanted form.

\vskip 0.3in

Assume first that ${\alpha}$ is a positive integer. Consider
$$H_n(Z_s,w)=H_n(Z_s,W_1)=\prod_{s\geq i\neq j\geq
1}(1-z_jz_i^{-1})^{n}\prod_{j=1}^s (1-z_jw^{-1})^{-n},$$ where we
identified $w_1$ with $w$ for simplicity. Notice that
$$H_{n}(Z_s,W_t)=H_{n}(Z_s,w)
\prod_{i=2}^t\prod_{j=1}^s(1-z_jw_i^{-1})^{-n}.$$ To prove Lemma
\ref{L:contraction} we need the following:

\begin{theorem}
For $n,s\in\mathbb{Z}$$_{>0}$,$s\geq 2$, there are polynomials
$f_{i,j}$ in $z_kz_l^{-1}$'s $(1\leq k\neq l\leq s)$ such that:

\begin{equation}
H_{n}(Z_s,w)=\sum_{i=1}^s \sum_{j=1}^{n}(1-z_iw^{-1})^{-j} f_{i,j}.
\end{equation}
Moreover for each $i$, $f_{i,j}$ is a polynomial in $z_i$. \vskip
0.3in
\end{theorem}

 Proof: To prove the existence of $f_{i,j}$'s, we will use induction on $s$.

 In the case of $s=2$ it is true by Lemma \ref{L:split}. Assume that it holds true for $s$, we have

 $$H_{n}(Z_{s+1},w)=H_{n}(Z_s,w)A_{s+1}=\sum_{i=1}^s
\sum_{j=1}^{n}(1-z_iw^{-1})^{-j} A_{s+1}f_{i,j}
$$$$=\sum_{i=1}^s
\sum_{j=1}^{n}(1-z_iz_{s+1}^{-1})^j(1-z_iw^{-1})^{-j}(1-z_{s+1}z_i^{-1})^{n}(1-z_{s+1}w^{-1})^{-n}
B_{s+1,i,j}f_{i,j},$$ where
$$A_{s+1}=(1-z_{s+1}w^{-1})^{-n}\prod_{l=1}^s(1-z_lz_{s+1}^{-1})^{n}(1-z_{s+1}z_l^{-1})^{n}$$
$$=(1-z_iz_{s+1}^{-1})^j(1-z_{s+1}z_i^{-1})^{n}(1-z_{s+1}w^{-1})^{-n}
B_{s+1,i,j}.$$ Notice that the term inside the sum can be split by
Lemma \ref{L:split}, while $B_{s+1,i,j}$ is a product of
$(1-z_kz_l^{-1})$'s, the existence follows.\vskip 0.1 in

As for the second part, note that $H_{n}(Z,w)$ is symmetric about
$z_1,\cdots,z_s$, we only need to prove  that
$f_{1,j},(j=1,\cdots,n)$ are polynomials of $z_1$. Multiplying two
sides of (4.1) by $(1-z_1w^{-1})^{n}\cdots(1-z_sw^{-1})^{n}$, we
have
\begin{align} \nonumber
&\prod_{1\leq i\neq j\leq s}(1-z_jz_i^{-1})^{n}\\
\label{E:power}
&=\sum_{i=1}^s\sum_{j=1}^{n}(1-z_1w^{-1})^{n}\cdots
(1-z_iw^{-1})^{n-j}\cdots (1-z_sw^{-1})^{n}f_{i,j}.
\end{align}

Using induction on $j'=n-j$: first, let $w=z_1$ in Eq.
(\ref{E:power}), we have

$$\prod_{1\leq i\neq j\leq
s}(1-z_iz_j^{-1})^{n}=f_{1,n}\prod_{i=2}^s(1-z_iz_1^{-1})^{n}$$

Eliminating the common factor we
find,$$f_{1,n}=\prod_{i=2}^s(1-z_1z_i^{-1})^{n}\prod_{2\leq i\neq
j\leq s}(1-z_iz_j^{-1})^{n}$$ which implies the case $j'=0$. Assume
that it's true for $j'<r$. Let $j'=r\leq n-1$. Differentiating both
sides of Eq. (\ref{E:power}) with respect to $z_1$, and set $w=z_1$,
we have:

$$\frac{\partial^r}{\partial z_1^r}\prod_{1\leq i\neq j\leq
s}(1-z_iz_j^{-1})^{n}=\prod_{i=2}^s(1-z_i/z_1)^{n}\cdot\sum_{i=0}^r
\binom{r}{i}i!(-z_1^{-1})^i\frac{\partial^{r-i}}{\partial
z_1^{r-i}}f_{1,n-i}.$$

The term $i=r$ in the sum contains $f_{1,n-r}$ and one finds that,
\begin{eqnarray*}
f_{1,n-r}&=&(r!)^{-1}(-z_1)^r\prod_{i=2}^s(1-z_i/z_1)^{-n}\frac{\partial^r}{\partial
z_1^r}\prod_{1\leq i\neq j\leq
s}(1-z_iz_j^{-1})^{n}\\
&-&\sum_{i=0}^{r-1}((r-i)!)^{-1}(-z_1)^{r-i}\frac{\partial^{r-i}}{\partial
z_1^{r-i}}f_{1,n-i}.
\end{eqnarray*}

Note that $$\frac{\partial^r}{\partial z_1^r}\prod_{1\leq i\neq
j\leq s}(1-z_iz_j^{-1})^{n}$$

$$=\prod_{2\leq i\neq j\leq
s}(1-z_iz_j^{-1})^{n}\sum c(a_i,b_i)
\prod_{i=2}^s\frac{\partial^{a_i}}{\partial
z_1^{a_i}}\Big((1-z_1/z_i)^{n}\Big)\frac{\partial^{b_i}}{\partial
z_1^{b_i}}\Big((1-z_i/z_1)^{n}\Big),$$ where the sum is over vectors
$(a_2\cdots a_s,b_2,\cdots,b_s)$ with nonnegative integer components
which sum up to $r$, And $c(a_i,b_i)=r!/(a_2!\cdots a_s!b_2!\cdots
b_s!)$ Now the first part of $f_{1,n-r}$ is
$$(r!)^{-1}(-1)^r z_1^{r-b_2-\cdots-b_s}\prod_{2\leq i\neq j\leq
s}(1-z_iz_j^{-1})^{n}\cdot$$$$\sum c(a_i,b_i)
\prod_{i=2}^s\frac{\partial^{a_i}}{\partial
z_1^{a_i}}\Big((1-z_1/z_i)^{n}\Big)\frac{\partial^{b_i}}{\partial
z_1^{b_i}}\Big((1-z_i/z_1)^{n}\Big)(1-z_iz_1^{-1})^{-n}z_1^{b_i}$$

By the following lemma and the assumption of induction,
$\lim_{z_1\rightarrow 0}f_{1,n-r}$ exists if $z_i\neq 0$. Observe
that $f_{1,n-r}$ is polynomial of $z_kz_l^{-1}$'s, it should be a
polynomial of $z_1$ as well. \vskip 0.3in

\begin{lemma} Let
$g_{s,m}(z)=z^s(1-a/z)^{-m}\frac{\partial^s}{\partial
z^s}(1-a/z)^m$,$a\neq 0$, then $\lim_{z\rightarrow 0}g_{s,m}(z)$
exists for $s\geq 0,0\leq s<m$.\end{lemma}

Proof: We use induction on $s$ again. The initial step is trivial.
Consider $s+1<m$,
$$g_{s+1,m}(z)=z^{s+1}(1-a/z)^{-m}\frac{\partial^s}{\partial
z^s}\Big(m(1-a/z)^{m-1}az^{-2}\Big)$$
$$=z^{s+1}(1-a/z)^{-m}am\sum_{i=0}^sd_iz^{-(2+s-i)}\frac{\partial^i}{\partial
z^i}(1-a/z)^{m-1}$$
$$=\sum_{i=0}^sc_ig_{i,m-1}(z)(z-a)^{-1}, $$ where $c_i=am(^s_i)(1+s-i)!(-1)^{s-i}=amd_i$,
the lemma follows.

\bigskip

\bibliographystyle{amsalpha}

\end{document}